\documentclass{amsart}

\newcommand{\A}{\mathcal A}
\newcommand{\sch}{\mathrm{ht}}
\newcommand{\com}{\mathrm{cx}}
\newcommand{\I}{\mathcal I}
\newcommand{\U}{\mathcal U}
\newcommand{\V}{\mathcal V}
\newcommand{\w}{\omega}

\newtheorem{theorem}{Theorem}
\newtheorem{lemma}{Lemma}
\theoremstyle{definition}
\newtheorem{remark}{Remark}

\title{Uniform Eberlein compactifications of metrizable spaces}
\author{Taras Banakh, Arkady Leiderman}
\address{Uniwersytet Humanistyczno-Przyrodniczy Jana Kochanowskiego w Kielcah, Poland and Ivan Franko National University of Lviv, Ukraine}
\email{t.o.banakh@gmail.com}
\address{Ben-Gurion University of Negev, Israel}
\email{arkady@math.bgu.ac.il}
\begin{document}

\begin{abstract} We prove that each metrizable space $X$ (of size $|X|\le\mathfrak c$) has a (first countable) uniform Eberlein compactification and each scattered metrizable space has a scattered hereditarily paracompact compactification. Each compact scattered hereditarily paracompact space is uniform Eberlein and belongs to the smallest class $\A$ of compact spaces, that contain the empty set, the singleton, and is closed under producing the Aleksandrov compactification of the topological sum of a family of compacta from the class $\A$.
\end{abstract}
\maketitle

Looking for compactifications of metrizable spaces it is natural to search them among compacta having as much properties of metrizable spaces as possible. In this respect, the class of Eberlein compacta fits quite well: Eberlein compacta are Frechet-Urysohn, they contain metrizable dense $G_\delta$-subsets, their weight coincides with their cellularity etc. In \cite{Ar} A.Arkhangelskii showed that each metrizable space has an Eberlein compactification. 

The class of Eberlein compacta contains a strictly smaller subclass consisting of uniform Eberlein compacta. Let us recall that a compact space $K$ is ({\em uniform}) {\em Eberlein} if $K$ is homeomorphic to a compact subspace of a (Hilbert) Banach space endowed with the weak topology, see \cite{BS}. 

In this paper we improve the mentioned compactification result of A.Arkhangelskii \cite{Ar} constructing uniform Eberlein compactification for each metrizable space.

\begin{theorem} Each metrizable space $X$ has a uniform Eberlein compactification.
\end{theorem}

\begin{proof} Instead of modifying the original construction of Arkhangelskii \cite{Ar} we shall present an alternative analytic proof. 

Given a metric space $X$, apply the classical Dowker result \cite[4.4.K]{En} to find an embedding $f:X\to S$ to the unit sphere $S$ of a Hilbert space $H$ of suitable density. It is well-known  that the norm and weak topologies coincide on $S$. So, $S$ can be considered as a subspace of the closed unit ball $B$ of $H$ endowed with the weak topology. In such a way we obtain an embedding of the metric space $X$ into the metrizable subspace $S$ of the weak unit ball $B$ of $H$, which is a uniform Eberlein compact. The closure of $X$ in $B$ will be the desired uniform Eberlein compactification of $X$.
\end{proof}

Modifying the original approach of Arkhangelskii, we can construct first countable uniform Eberlein compactifications.

\begin{theorem} A metrizable space $X$ has a first countable uniform Eberlein compactification if and only if $|X|\le\mathfrak c$.
\end{theorem}

\begin{proof} The ``only if'' part follows from the famous Arkhangelskii's result \cite[3.1.30]{En} saying that each first countable compact Hausdorff space $K$ has cardinality $|K|\le\mathfrak c$. To prove the ``if' part, we shall use Kowalsky's Theorem \cite[4.4.9]{En} asserting that each metrizable space $X$ of weight $\kappa$ embeds into the countable power $(H_\kappa)^\w$ of the metric hedgehog
$$H_\kappa=\{t\vec e_\alpha:t\in[0,1],\;\alpha\in\kappa\}\subset l_2(\kappa)$$with $\kappa$ spines. Here $(\vec e_\alpha)_{\alpha\in\kappa}$ is the standard orthonormal base of the Hilbert space $l_2(\kappa)$. If $X$ has weight $\le\mathfrak c$, then $X$ embeds into the countable power $H_{\mathfrak c}^\omega$ of $H_{\mathfrak c}$. Since the countable product of first countable uniform Eberlein compacta is first countable and uniform Eberlein, it suffices to show that the hedgehog $H_{\mathfrak c}$ has a first countable uniform Eberlein compactification. This compactification can be constructed explicitly. First fix an injective enumeration $\{\varphi_\alpha:\alpha\in\mathfrak c\}$ of a topological copy of the Cantor set in $[\pi/6,\pi/3]$.

The embedding $f:H_{\mathfrak c}\to S$ of the hedgehog $H_{\mathfrak c}$ into the unit sphere $S$ of the Hilbert space $l_2(\mathfrak c+2)$ can be defined as follows:
$$f:t\vec e_\alpha\mapsto \cos(t)\vec e_{\mathfrak c+1}+\sin(t)\,\big(\cos(\varphi_\alpha)\vec e_\mathfrak c+\sin(\varphi_\alpha)
\vec e_\alpha\big)\quad \mbox{for $t\in[0,1]$ and $\alpha\in\mathfrak c$}.$$
 One can show that the closure of $f(H_{\mathfrak c})$ in the weak topology of $l_2(\mathfrak c+2)$ coincides with the set $$f(H_{\mathfrak c})\cup\{\cos(t)\vec e_{\mathfrak c+1}+\sin(t)\cos(\varphi_\alpha)\vec e_\mathfrak c:t\in[0,1],\;\alpha\in\mathfrak c\}$$
homeomorphic to the cone over the Aleksandrov duplicate of the Cantor set.
The latter set is first countable and being a compact subset  of $(l_2(\mathfrak c+2),\mathit{weak})$, is uniform Eberlein.
\end{proof}

\begin{remark} Following the approach of \cite{KM} or \cite{C} one can prove that each metrizable space $X$ (of cardinality $|X|\le\mathfrak c$) has a (first countable) uniform Eberlein compactification $K$ with $\dim(K)=\dim(X)$.
\end{remark}

Next, we shall be interested in scattered uniform Eberlein compactifications of scattered metrizable spaces. We recall that a topological space $X$ is  {\em scattered} if each subspace of $X$ has an isolated point. According to \cite{Tel} each metrizable scattered space has a compactification homeomorphic to an initial segment of ordinals endowed with the order topology. Unfortunately, such uncountable segments are not Fr\'echet-Urysohn and hence not Eberlein compact. Nonetheless,  scattered metrizable spaces have scattered compactifications which are hereditarily paracompact and uniform Eberlein. 

The class of scattered hereditarily paracompact compacta is very narrow and admits a constructive description. Namely, it coincides with the smallest class $\A$ of compacta, which contains the empty set, the singleton, and the Aleksandrov compactifications of topological sums of compacta from the class $\A$. We recall that a compactification $\overline X$ of a Tychonov locally compact space $X$  is {\em Aleksandrov} if its remainder  has size $|\overline{X}\setminus X|\le 1$. 

\begin{theorem}\label{main}\begin{enumerate}
\item Each scattered metrizable space $X$ has a scattered hereditarily paracompact compactification $K$.
\item Each scattered hereditarily paracompact compact space is uniform Eberlein.
\item  A compact space $X$ is scattered and hereditarily paracompact if and only if it belongs to the class $\A$. 
\end{enumerate}
\end{theorem}

\begin{remark} The class $\A$ is not closed with respect to finite products: the product $\alpha\aleph_1\times\alpha\aleph_0$ of the Aleksandrov compactifications of the discrete spaces of size $\aleph_1$ and $\aleph_0$ does not belong to the class $\A$ because it is not hereditarily normal, see \cite[2.7.16(a)]{En}. On the other hand, $\alpha\aleph_1\times\alpha\aleph_0$ is a scattered uniform Eberlein compact space. It is known that each scattered Eberlein compact space $K$ is strong Eberlein in the sense that it embeds into a $\sigma$-product $$\{(x_\alpha)_{\alpha\in\kappa}\in\{0,1\}^\kappa:|\{\alpha\in\kappa:x_\alpha\ne0\}|<
\aleph_0\}$$
of the two-point space for some cardinal $\kappa$, see \cite{Al}. 
\end{remark}

\begin{remark} In light of Theorem~\ref{main} it should be mentioned that there are scattered countable (and hence hereditarily paracompact) Tychonov spaces admitting no scattered compactification, see \cite{Mal}, \cite{Nyi}, \cite{Ra}, \cite{Sol}.
\end{remark}

The remaining part of the paper is devoted to the proof of Theorem~\ref{main}, which is done by transfinite induction using two ordinal functions: the scattered height of a scattered space and the complexity of a space $K\in\mathcal A$.

The scattered height of a scattered space $X$ is defined using the $\alpha$-th derived sets $X^{(\alpha)}=\bigcap_{\beta<\alpha}(X^{(\beta)})'$ where $X^{(0)}=X$ and $A'$ denotes the set of non-isolated points of a space $A$. It is easy to see that a space $X$ is scattered if and only if $X^{(\alpha)}$ is empty for some ordinal $\alpha$. The smallest ordinal $\alpha$ such that $|X^{(\alpha)}|\le 1$ is called the {\em scattered height} of $X$ and is denoted by $\sch(X)$. For example, a space $X$ with a unique non-isolated point has $\sch(X)=1$.

Next, we observe that the class $\A$ can be written as the union $\mathcal A=\bigcup_\alpha\mathcal A_\alpha$ where $\mathcal A_0$ is the class of all spaces $X$ with $|X|\le 1$, and $\mathcal A_\alpha$ consists of the Aleksandrov compactifications of topological sums of compacta from the class $\A_{<\alpha}=\bigcup_{\beta<\alpha}\mathcal A_\beta$. The {\em complexity} $\com(K)$ of a space $K\in\mathcal A$ is the smallest ordinal $\alpha$ such that $K\in\mathcal A_\alpha$. 

\begin{lemma}\label{l1} Any compact space $K\in\mathcal A$ is a scattered, hereditarily paracompact, uniform Eberlein and has $\com(K)=\sch(K)$.
\end{lemma}

\begin{proof} This lemma will be proven by induction on $\com(K)$. The assertion is trivial if $\com(K)=0$ (in which case $|K|\le1$).

Assume that the lemma is proved for compacta $X\in\A$ with $\com(X)<\alpha$.
Fix any compact space $K\in\A$ with $\com(K)=\alpha>0$. Then there is a family $\{K_i:i\in\I\}\subset\A_{<\alpha}$ such that $K$ is the Aleksandrov compactification of the topological sum $\oplus_{i\in\I}K_i$. By the inductive assumption, each space $K_i$  is scattered, hereditarily paracompact, uniform Eberlein and has $\com(K_i)=\sch(K_i)$. Then 
$$\com(K)=\sup_{i\in\I}(\com(K_i)+1)=\sup_{i\in\I}(\sch(K_i)+1)=\sch(K).$$

To show that $K$ is scattered and hereditarily paracompact, take any infinite subspace $X\subset K$. Being infinite, the set $X$ meets some set $K_i$. Since $K_i$ is scattered, the intersection $X\cap K_i$ contains an isolated point which is isolated in $X$ too (because $X\cap K_i$ is  open in $X$). This proves that $K$ is scattered. 

Next, we show that $X$ is paracompact. This is trivial if $X\subset\oplus_{i\in\I}K_i$ (because the hereditary paracompactness is preserved by topological sums). So we assume that $X$ contains the compactifying point $\infty$ of $K$. Take any open cover $\U$ of $X$ and find a set $U_\infty\in\U$ containing the compactifying point $\infty$. By the regularity of $X$, find a closed neighborhood $W\subset U_\infty$ of $\infty$. The hereditary paracompactness of $\oplus_{i\in\I}K_i$ allows us to find a locally finite open cover $\V$ of $X\setminus\{\infty\}$ refining the cover $\U$. Then $\mathcal W=\{U_\infty\}\cup\{V\setminus W:V\in\V\}$ is a locally finite open cover of $X$ refining the cover $\U$ and witnessing that $X$ is paracompact.

Finally, we show that $K$ is a uniform Eberlein compact. For every $i\in\I$ find an embedding $e_i:K_i\to H_i$ into an infinite-dimensional Hilbert space endowed with the weak topology. Take a family $\{\Gamma_i\}_{i\in\I}$ of pairwise disjoint sets of cardinality $|\Gamma_i|$ equal to the density of the Hilbert space $H_i$. By the classical Riesz theorem, $H_i$ can be identified with the Hilbert space $l_2(\Gamma_i)=\{f:\Gamma_i\to\mathbb{R}:\sum_{\gamma\in\Gamma_i}f(\gamma)^2<+\infty\}$.
After a suitable affine transformation, we can assume that the image $e_i(K_i)$ does not contain the origin $\mathbf{0}$ of $l_2(\Gamma_i)=H_i$ and lies in the closed unit ball centered at the origin. Let $\Gamma=\bigcup_{i\in\I}\Gamma_i$ and observe that the map $e:K\to l_2(\Gamma)$ defined by
$$e(x)=\begin{cases}e_i(x)&\mbox{if $x\in K_i$ for some $i\in\I$}\\
\mathbf{0}&\mbox{if $x\in K\setminus\oplus_{i\in\I}K_i$}
\end{cases}
$$is an embedding of $K$ into the Hilbert space $l_2(\Gamma)$ endowed with the weak topology. This witnesses that  $K$ is a uniform Eberlein compact.
\end{proof}

A topological space $X$ is called {\em strongly zero-dimensional} if each open cover $\U$ of $X$ has a discrete open refinement $\V$ (the latter means that each point of $X$ has a neighborhood meeting at most one set from $\V$). We shall need the following result due to R.~Telgarsky \cite{Tel}.

\begin{lemma}\label{l2} Each scattered paracompact space is strongly zero-dimensional.
\end{lemma}

In the following lemma we prove the third statement of Theorem~\ref{main}. This lemma combined with Lemma~\ref{l1} implies also the second statement  of Theorem~\ref{main}.

\begin{lemma}\label{l3} A compact space $X$ belongs to the class $\A$ if and only if $X$ is scattered and hereditarily paracompact.
\end{lemma}

\begin{proof} The ``only if'' part has been proved in Lemma~\ref{l1}. To prove the ``if'' part, assume that $X$ is a compact, scattered, hereditarily paracompact space.
We need to prove that $X\in\A$. This will be done by induction on the scattered height of $X$. The inclusion $X\in\A$ is trivial if $\sch(X)=0$ (in which case $X$ is empty or a singleton). Assume that $\sch(X)=\alpha>0$ and all scattered hereditarily paracompact spaces with scattered height $<\alpha$ belong to the class $\A$.

It follows that the $\alpha$-derived set $X^{(\alpha)}$ contains at most one point. If it is empty, then $\alpha$ is not limit and $X^{(\alpha-1)}$ is finite. Being scattered, the compact space $X$ is zero-dimensional. Consequently, we can find a finite cover $\{X_i:i\in\I\}$ of $X$ consisting of pairwise disjoint clopen subsets such that each set $X_i$ has at most one-point intersection with the set $X^{(\alpha-1)}$. It follows that $|X_i^{(\alpha-1)}|\le 1$ and thus $\sch(X_i)\le\beta-1<\alpha$. By the inductive assumption, each $X_i\in\A$ and then $X\in\A$, being the Aleksandrov compactification of the topological sum $\oplus_{i\in\I}X_i=X$.

Next, we consider the case of a singleton $X^{(\alpha)}=\{a\}$.  In this case $X$ is the Aleksandrov compactification of the locally compact space $X\setminus\{a\}$. The space $X\setminus\{a\}$, being scattered and paracompact, is strongly zero-dimensional by Lemma~\ref{l2}. This allows us to find a  discrete cover $\{X_i:i\in\I\}$ of $X\setminus\{a\}$ that consists of clopen subsets of $X\setminus\{a\}$ whose closures in $X$ do not contain the compactifying point $\{a\}$. Moreover, if the ordinal $\alpha$ is not limit, we may additionally require that each set $X_i$, $i\in\I$ has at most one-point intersection with the set $X^{(\alpha-1)}$ (which has a unique non-isolated point $a$ and hence is closed and discrete in $X\setminus\{a\}$).
It follows that each space $X_i$ is compact and has $\sch(X_i)<\alpha$. By the inductive assumption, each space $X_i$ belong to $\A$. The discreteness of the family $\{X_i\}$ in $X\setminus\{a\}$ implies that $X=\{a\}\cup\bigcup_{i\in\I}X_i$ can be identified with the Aleksandrov compactification of the topological sum $\oplus_{i\in\I}X_i$. Consequently, $X\in\A$ by the definition of $\A$.
\end{proof}

It is well-known that each ordinal $\alpha$ can be uniquely written as $\alpha=\beta+n(\alpha)$ for some limit ordinal $\beta$ and some finite ordinal $n(\alpha)$.

The following lemma combined with Lemma~\ref{l2} yields (a controlled version of) the first statement of Theorem~\ref{main}.

\begin{lemma}\label{l4} Each  scattered metrizable space $X$ with scattered height $\sch(X)=\alpha$ has a compactification $K\in\A$ with scattered height $\sch(K)\le \alpha+n(\alpha)+1$.
\end{lemma}

\begin{proof} The proof will be done by induction on $\alpha$. If $\sch(X)=0$, then $|X|\le 1$ and $K=X$ is the compactification of $X$ in the class $\A$ with $\sch(K)=0\le 1=\alpha+n(\alpha)+1$. 

Now assume that the lemma is proved for scattered metrizable spaces with scattered height $<\alpha$. Let $X$ be a scattered metrizable space with $\sch(X)=\alpha$. Then $|X^{(\alpha)}|\le1$ but $|X^{(\beta)}|>1$ for all $\beta<\alpha$.
 
We shall consider separately three cases.
\smallskip

1. $\alpha$ is a limit ordinal and $X^{(\alpha)}=\emptyset$. Then $\bigcap_{\beta<\alpha}X^{(\beta)}=X^{(\alpha)}=\emptyset$ and hence $X=\bigcup_{\beta<\alpha}X\setminus X^{(\beta)}$. Being scattered and paracompact, the space $X$ is strongly zero-dimensional by Lemma~\ref{l2}. Consequently, it admits a discrete open cover $\{U_i:i\in\I\}$ such that each set $U_i$ lies in $X\setminus X^{(\beta)}$ for some $\beta<\alpha$. Then $U_i^{(\beta)}=\emptyset$ and hence $\sch(U_i)\le\beta<\alpha$ for all $i\in\I$. By the inductive assumption, $U_i$ has a compactification $K_i\in\A$ with $\sch(K_i)\le\beta+n(\beta)+1<\alpha$. Let $K$ be the Aleksandrov compactification of the topological sum $\oplus_{i\in\I}K_i$. It follows  that $\sch(K)\le \alpha=\alpha+n(\alpha)$. 
The discreteness of the cover $\{U_i\}_{i\in\I}$ implies that $X$ is homeomorphic to the topological sum $$\oplus_{i\in\I}U_i\subset \oplus_{i\in \I}K_i\subset K$$ and hence $K$ is the compactification of $X$.
\smallskip

2. $\alpha$ is not limit and $X^{(\alpha)}=\emptyset$. Then $X^{(\alpha-1)}$ is a closed discrete subspace of $X$. Since $X$ is strongly zero-dimensional, there is a discrete cover $\{X_i\}_{i\in\w}$ of $X$ by clopen subsets having at most one-point intersection with the set $X^{(\alpha-1)}$. It follows that for each $i\in\I$ we get $|X_i^{(\alpha-1)}|\le 1$ and thus $\sch(X_i)\le\alpha-1<\alpha$. By the inductive assumption, the space $X_i$ has a compactification $K_i\in\A$ with $\sch(K_i)\le\alpha-1+n(\alpha-1)+1=\alpha+n(\alpha)-1$. The Aleksandrov compactification $K$ of the topological sum $\oplus_{i\in\I}K_i$ belongs to $\A$ and has scattered height $\sch(K)\le \alpha+n(\alpha)$. Since $X=\oplus_{i\in\I}X_i$, $K$ is a compactification of $X$.
\smallskip

3. $X^{(\alpha)}=\{a\}$ is a singleton. Being first countable, $X$ has a decreasing neighborhood base $(W_n)_{n\in\w}$ at $x$ with $W_0=X$. Moreover, taking into account that $X$ is zero-dimensional (being scattered and paracompact), we may assume that each set $W_n$ is clopen in $X$. For every $n\in\w$ consider the clopen subspace $Y_n=W_n\setminus W_{n+1}$ of $X$. Since $Y_n\cap X^{(\alpha)}=\emptyset$, we conclude that $Y_n^{(\alpha)}=\emptyset$ and hence $\beta=\sch(Y_n)\le\alpha$. If $\beta=\alpha$, then by the preceding cases we can find a compactification $K_n\in\A$ of $Y_n$ with $\sch(K_n)\le \alpha+n(\alpha)$. If $\beta<\alpha$, then by the inductive assumption we can find a compactification $K_n\in\A$ of $Y_n$ with $\sch(K_n)\le\beta+n(\beta)+1\le\alpha+n(\alpha)$. The Aleksandrov compactification $K$ of the topological sum $\oplus_{n\in\w}K_n$ belongs to $\A$, is a compactification of the space $X=\{a\}\bigcup_{n\in\w}Y_n$, and has scattered height $\sch(K)\le \sup_{n\in\w}(\sch(K_n)+1)\le\alpha+n(\alpha)+1$.
\end{proof}


\end{document}